\documentclass[10pt]{amsart} 

\usepackage{amsfonts,amssymb,amsthm,amsmath}
\usepackage{eucal}
\usepackage{ifthen,epsfig,color,graphics,rotating,hangcaption}

\newcommand{\showcomments}{yes}

\newsavebox{\commentbox}
{\ifthenelse{\equal{\showcomments}{yes}}%
{\footnotemark
    \begin{lrbox}{\commentbox}
    \begin{minipage}[t]{1in}\raggedright\sffamily\small
    \footnotemark[\arabic{footnote}]}
{\begin{lrbox}{\commentbox}}}%
{\ifthenelse{\equal{\showcomments}{yes}}%
{\end{minipage}\end{lrbox}\marginpar{\usebox{\commentbox}}}
{\end{lrbox}}}

\newtheorem{thm}{Theorem}[section]
\newtheorem{prop}[thm]{Proposition}
\newtheorem{cor}[thm]{Corollary}

\theoremstyle{definition}
\newtheorem{defn}[thm]{Definition}
\newtheorem{exmp}[thm]{Example}

\newtheorem{que}[thm]{Question}
\newtheorem{prob}[thm]{Problem}

\theoremstyle{remark}

\newtheorem{rem}{Remark}

\renewcommand{\Re}{\text{Re}}

\renewcommand{\(}{\left(}
\renewcommand{\)}{\right)}  
\newcommand{\abs}[1]{\left\arrowvert #1\right\arrowvert}
\newcommand{\norm}[1]{\left\| #1\right\|}

\newcommand{\supp}{\text{supp}}

\newcommand{\field}[1]{\mathbb{#1}}
 
\newcommand{\FF}{\ensuremath{\field{F}}}
\newcommand{\KK}{\ensuremath{\field{K}}} 
\newcommand{\II}{\ensuremath{\field{I}}} 
\newcommand{\IF}{\ensuremath{\field{I}\field{F}}} 
\newcommand{\IR}{\ensuremath{\field{I}\field{R}}} 
 
\newcommand{\CC}{\ensuremath{\field{C}}} 
\newcommand{\Ct}{\ensuremath{\field{C}^2}}
\newcommand{\Cn}{\ensuremath{\field{C}^n}}
\newcommand{\RR}{\ensuremath{\field{R}}}
\newcommand{\Rt}{\ensuremath{\field{R}^2}}
\newcommand{\Rn}{\ensuremath{\field{R}^n}}

\newcommand{\eps}{\epsilon}

\newcommand{\boxcov}{\mathcal{B}}
\newcommand{\chainrec}{\mathcal{R}}

\newcommand{\chreceps}{\chainrec_{\eps}}
\newcommand{\chrecepsp}{\chainrec_{\eps'}}

\newcommand{\boxchrecmod}{box chain recurrent model } 

\newcommand{\drawfigzcmplxcant}{\scalebox{1}{\includegraphics{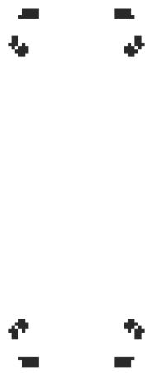}}}
\newcommand{\drawfigcantcirc}{\scalebox{1}{\includegraphics{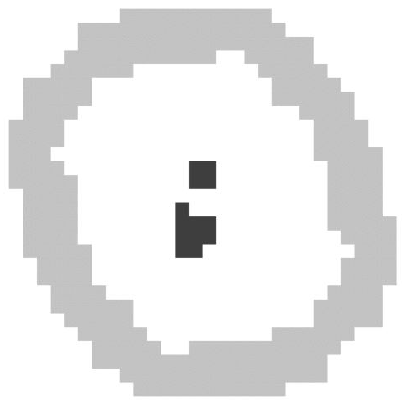}}}

\newcommand{\drawfigcantxbasb}{\scalebox{1}{\includegraphics{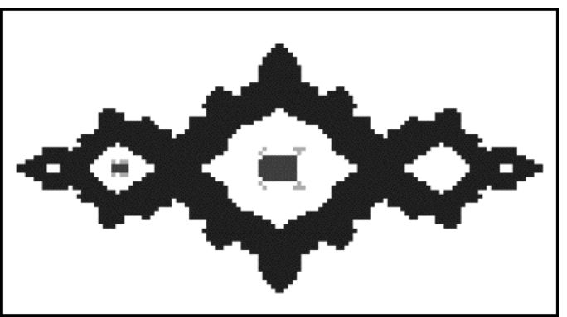}}}
\newcommand{\drawfigcantxbasnegb}{\scalebox{1}{\includegraphics{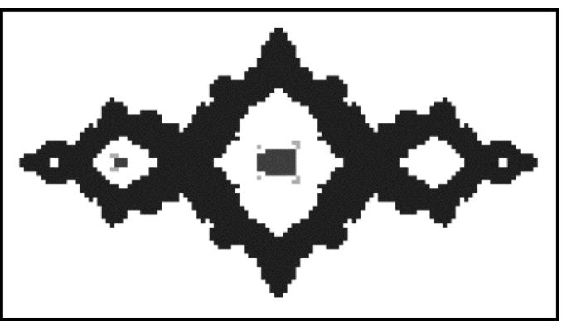}}}

\newcommand{\drawfigmatcircnega}{\scalebox{1}{\includegraphics{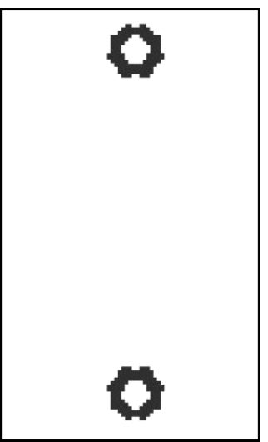}}}
\newcommand{\drawfigmatcircb}{\scalebox{1}{\includegraphics{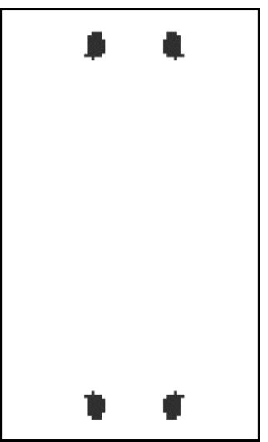}}}
\newcommand{\drawfigmatcirca}{\scalebox{1}{\includegraphics{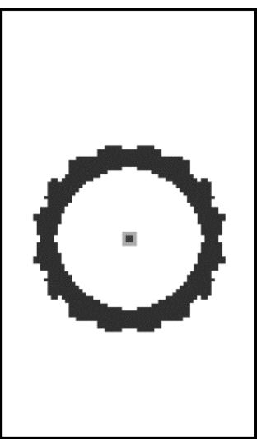}}}
\newcommand{\drawfigmatcircnegb}{\scalebox{1}{\includegraphics{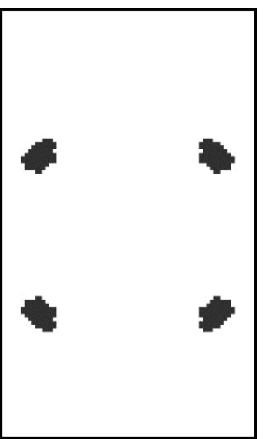}}}

\newcommand{\drawfigmatrabnega}{\scalebox{1}{\includegraphics{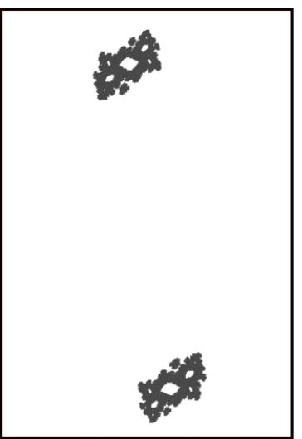}}}
\newcommand{\drawfigmatrabb}{\scalebox{1}{\includegraphics{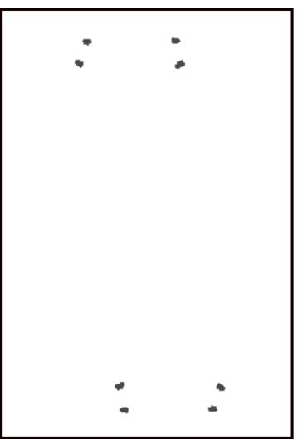}}}
\newcommand{\drawfigmatraba}{\scalebox{1}{\includegraphics{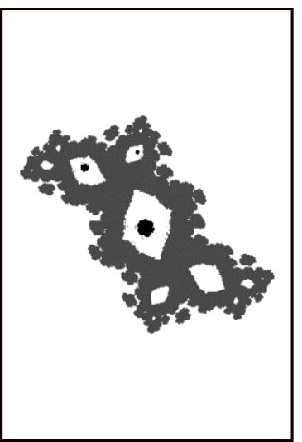}}}
\newcommand{\drawfigmatrabnegb}{\scalebox{1}{\includegraphics{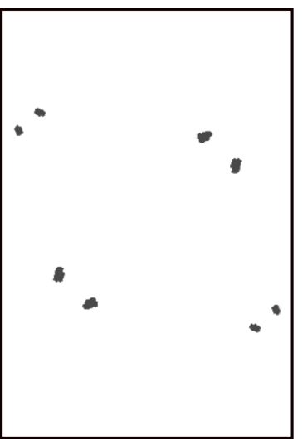}}}

\title[Modeling polynomial skew products]
{Rigorous numerical studies of the \\ dynamics of
 polynomial skew products of $\Ct$
 }
\author[S.L. ~Hruska]{Suzanne Lynch Hruska}
\thanks{Research supported in part by a grant from the National Science Foundation.}
\address{Department of Mathematics\\
Indiana University\\
Rawles Hall\\
Bloomington, IN 47405, USA
}
\email{shruska@msm.umr.edu}


\date{\today}

\begin{document}

\begin{abstract} For the class of polynomial skew products of $\Ct$, we 
describe a rigorous computer algorithm which, for a given map $f$, will 
(1) build a model of the dynamics of $f$ on its
chain recurrent set, and (2) attempt to determine whether
$f$ is Axiom A.  Further, we discuss how we used our implementation of this algorithm
 to establish Axiom A for several explicit cases.
\end{abstract}

\maketitle

\section{Introduction}
\label{sec:intro}

Our main interest here is to develop and use rigorous computer investigations to study the dynamics of polynomial skew products of $\Ct$; i.e., maps of the form
$
f(z,w) = (p(z), q(z,w)),
$
where $p$ and $q$ are polynomials of the same degree $d\geq 2$.

The skew products we are most interested in studying 
are those maps which are {\em Axiom A}.  Such maps have the ``simplest'' chaotic dynamics, and stability under small perturbation, thus are amenable to computer investigation.  
In one complex dimension, a polynomial map is called {\em hyperbolic} if it is uniformly expanding on some neighborhood of its Julia set, with respect to some riemannian metric.  In $\Cn$, a map is {\em hyperbolic on an invariant set $\Lambda$} if there exists a continuous splitting of the tangent bundle over $\Lambda$ into two subspaces (of any dimension zero through $n$), with one subspace uniformly expanded by the map, and the other uniformly contracted. A map is {\em Axiom A} if $f$ is hyperbolic on the nonwandering set, $\Omega$, and if periodic points are dense in $\Omega$ and $\Omega$ is compact.

Our motivation is to understand what kind of dynamics can occur for Axiom A polynomial skew products.  
In this paper, we construct a class of maps with interesting dynamics, and prove using rigorous computer techniques that sample maps from this class are Axiom A.  This leads us to conjecture that all (or nearly all) maps in this class are Axiom A.

To develop our computer techniques, we use the work of \cite{SLHone} and \cite{SLHtwo} as a foundation.  In \cite{SLHone},
 we described a 
rigorous algorithm (and its implementation) for constructing a neighborhood $\boxcov$ of 
the chain recurrent set, $\chainrec$,
and a graph $\Gamma$ modelling the dynamics of $f$ on $\boxcov$.
In \cite{SLHtwo} we developed a numerical method for proving hyperbolicity of a polynomial map $p$ of $\CC$, which consists of  an attempt to construct a metric in which $p$ is expanding on $\boxcov$, by some uniform factor $L>1$.

Here, we adapt the one-dimensional algorithm to the setting of skew products.  Our efforts are aided by the fact that skew products are a natural generalization of polynomials in one dimension, since a skew product $f$ maps the vertical line $\{z\} \times \CC$ to the vertical line $\{p(z)\} \times \CC$, and restricted to a vertical line, $f$ is the polynomial map $w \mapsto q_z(w)=q(z,w)$.  

We applied a similar approach in \cite{SLHthree}, where we developed a test for verifying hyperbolicity of polynomial diffeomorphisms of $\Ct$.  This algorithm is different from that for skew products, since the dynamics of polynomial diffeomorphisms is of saddle type, i.e., with one direction expanding and the other contracting, whereas for skew products, it happens that to establish Axiom A we need only check an expanding condition.  

We have implemented all of our algorithms  into a computer program called {\em Hypatia}, which began with the algorithms of \cite{SLHone}, then was enhanced to include algorithms from \cite{SLHtwo} (and \cite{SLHthree}), and now with the present work, can be used to prove Axiom A for specific polynomial skew products.

Finally, we describe the organization of this paper. 
Polynomial skew products of $\Ct$ have been studied by Jonsson (\cite{MatSk, MatDil}) and Heinemann (\cite{HeinETDS1, HeinETDS2, HeinKy}).   
 In Section~\ref{sec:bckgnd}, we provide necessary background on skew products, hyperbolicity, and invariant sets of interest.  In Section~\ref{sec:newclasses}, we describe the dynamics of several classes of skew products, including our new class discussed above.
 In Section~\ref{sec:hypatiatheory}, we outline our rigorous computer algorithms, implemented in {\em Hypatia}, for attempting to verify Axiom A for skew products.  
In Section~\ref{sec:examples}, we describe how we  used {\em Hypatia} to prove Axiom A for several specific skew products, including maps of the same type as Jonsson's and Heinemann's examples, in addition to our new example.
  In Appendix~\ref{sec:IA},  we give a brief overview of the technique that we used to control round-off error in our computations, called {\em interval arithmetic}.

\subsection*{Acknowledgements}  
I would like to thank John Smillie for suggesting this project, Mattias Jonsson and Manfred 
Denker for advice on getting started, and Eric Bedford for guidance throughout my
investigations.  I would also 
like to thank Adrien Douady and Mikhail Lyubich for asking enlightening questions
when I presented this work at the Snowbird meeting in June 2004, and Amie Wilkinson
for mentioning the relationship between polynomial skew products and partial hyperbolicity.  

\section{Background}
\label{sec:bckgnd}


\subsection{Chain recurrence and hyperbolicity}
\label{sec:dpspaces}

The {\em chain recurrent set}, $\chainrec$,
the {\em Julia set}, $J$, and the {\em non-wandering set}, $\Omega$,  
are all attempts at
locating the points with dynamically interesting behavior.
The Julia set, $J$, can be defined as the topological boundary of the set, $K$, of points in $\Ct$ with bounded orbits under $f$.  Both $J$ and $K$ are compact.  
Slices of the Julia set can be easily sketched by computer, but
$\chainrec$ is the set most amenable to
rigorous computer investigation. $\chainrec$ can also be easily decomposed into components which do not interact with one another. 
Since $\chainrec \supset J$, 
we can learn about $J$ by studying $f$ on $\chainrec$.

%
An {\em $\eps$-chain} of length $n>1$ from $y$ to $z$ is a 
  sequence of points
  $\{y=x_1, \ldots , x_n=z\}$ such that  $\abs{f(x_k) - 
  x_{k+1}} < \epsilon$ for $1 \leq i \leq n-1.$
A point $y$ belongs to the {\em $\eps$-chain recurrent set}, 
  $\chreceps$, of a function~$f$ if there is an
  $\eps$-chain from $y$ to $y$.
The {\em chain recurrent set} is $\chainrec = \cap_{\eps>0} 
\chreceps.$
A point $z$ is in the {\em forward chain limit set of a point $y$},
$\chainrec(y)$, if for all $\eps >0$, for all 
$n\geq 1$, there is an $\eps$-chain from $y$ to $z$ of length 
greater than $n$.
Put an equivalence relation on $\chainrec$ by: 
$y \sim z$ if $y \in \chainrec(z)$ and $z\in\chainrec(y)$.
Equivalence classes are called {\em chain transitive components}.
Define $\chreceps(y)$ and $\eps$-chain transitive
components analogously. 
 $\chainrec$ is closed and invariant, and if $\eps < \eps'$,
then $\chainrec \subset \chreceps
 \subset \chrecepsp$.  

Thus chain recurrence is quite natural to rigorously study using a computer. 
Next, we recall the precise definitions of hyperbolicity and Axiom A.

\begin{defn} \label{defn:hyp}
Let $g$ be a $C^1$ diffeomorphism or endomorphism of a compact manifold $M$, and let $\Lambda$ be a closed, $g$-invariant set. 
Say $\Lambda$
is {\em hyperbolic} for $g$ if there is a splitting of the tangent bundle $T_xM
= E^s_x \oplus E^u_x$ (one subspace may be trivial), for each $x$ in $\Lambda$, which varies
continuously with $x$ in $\Lambda$, a constant $\lambda > 1$, and a
riemannian norm $\norm{\cdot}$ such that:
\begin{enumerate}
 \item $g$ preserves the splitting, \textit{i.e.}, $D_xg(E^s_x) = 
E^s_{gx}$, and $D_xg(E^u_x) = E^u_{gx}$, and 
\item $Dg$ expands (contracts) $E^u (E^s)$ uniformly, \textit{i.e.},
$\norm{D_xg(\mathbf{w})} \geq \lambda 
\norm{\mathbf{w}}$ for $\mathbf{w} \in E^u_x$, and 
 $\norm{D_xg(\mathbf{v})} \leq (1/\lambda) \norm{\mathbf{v}}$ for $\mathbf{v} \in E^s_x$.
\end{enumerate}
\end{defn}

This definition is independent of choice of norm.  

\begin{defn} \label{defn:axiomA}
$g$ is {\em Axiom A} if for the nonwandering set $\Omega$, we have: (1) $\Omega$ is compact, (2) periodic points are dense in $\Omega$, and (3) $\Omega$ is hyperbolic for $g$.
\end{defn}

If Per$(g)$ is the set of periodic points of $g$, then $\overline{\text{Per}(g)} \subset \Omega \subset \chainrec$.   Jonsson shows that for Axiom A polynomial skew products, $\overline{\text{Per}(f)} = \Omega = \chainrec$ (see Theorem~\ref{thm:AxiomAVertExp}).

\subsection{Polynomial skew products of $\CC^2$}
\label{sec:skewintro}

In this subsection, we summarize some of the notation and results of \cite{MatSk}, to give needed background on skew products.

Since the dynamics of $f$ in the $z$-coordinate is given by $p$, it will be useful to employ the notation ($K_p$ and) $J_p$ for the one-dimensional (filled) Julia set of $p$, and $G_p(z)$ for the Green function in $\CC$ of $p$, where $K_p = \{ G_p=0\}$.  

Let $\mathcal{M}$ denote the Mandelbrot set for quadratic polynomials, i.e., $\mathcal{M}$ is the set of all $c$ such that the critical orbit is bounded under $p_c(x) = x^2+c$ (equivalently, $\mathcal{M}$ is the locus of $c$ with connected $J_c$).   

%

\textbf{Global Dynamics.} For polynomial skew products,  the usual rate of escape Green function, defined for $x\in\Ct$ by $G(x)=\lim_{n\to\infty} \frac{1}{d^n} \log^+ \abs{f(x)}$,  is continuous, plurisubharmonic, nonnegative, and satisfies $G\circ f = dG$ and  $K=\{G=0\}$.
One can also define a positive closed current $T = \frac{1}{2\pi}dd^c G$ and an ergodic invariant measure, $\mu = T \wedge T$, of maximal entropy $\log d^2$.  
%
 $J_2:=\supp(\mu)$ is also the closure of the set of repelling periodic points.  
 %

\textbf{Vertical dynamics.}  Since $f$ preserves the vertical lines $\{z\} \times \CC$, it is useful to consider the dynamics of $f$ on this family of lines.  Let $z_n = p^n(z)$, $q_z(w) = q(z,w)$,  and $Q^n_z(w) = q_{z_{n-1}} \circ \cdots \circ q_z$, so that $f^n(z,w) = (z_n, Q^n_z(w)).$  Let $G_z(w) = G(z,w) - G_p(z)$.  Then $G_z$ is nonnegative, continuous, subharmonic, and is asymptotic to $\log\abs{w} - G_p(z)$ as $w\to\infty$.  
Naturally, define $K_z = \{ G_z=0\}$, and $J_z=\partial K_z$.  Then $K_z$ and $J_z$ are compact, and  if $z\in K_p$, then $w\in K_z$ if and only if $\abs{Q^n_z(w)}$ is bounded.  Further, $G_{z_1} \circ q_z = d G_z$, which implies $q_z(K_z) = K_{z_1}$ and  $q_z(J_z) = J_{z_1}$.
Define $C_z$ as the critical set in $\{z\} \times \CC$. For example, if $q(z,w) = w^2+b+\eps z$, then $C_z = \{0\}$ for every $z$.   For all $z\in\CC$, $J_z$ is connected if and only if $C_{z_n} \subset K_{z_n}$.

However, not every phenomena of one-dimensional dynamics carries over to vertical dynamics.  For example, unlike in one dimension, $J_z$ may have finitely many (but greater than one) connected components, even for $d=2$ (see \cite{MatSk}, remark 2.5).

%

\textbf{Vertical Expansion.}  Let $Z \subset K_p$ be compact with $pZ \subset Z$, for example $Z=J_p$ or $Z=A_p$, the set of attracting periodic orbits.  Let $J_Z = \overline{\cup_{z\in Z}\{z\} \times J_z }$.  Jonsson shows $J_{J_p} = J_2$.
Call $f$  {\em vertically expanding over $Z$} if there exist $c>0$ and $\lambda >1$ such that $\abs{DQ^n_z(w)} \geq c\lambda^n$, for all $z\in Z$, $w\in J_z$, and $n\geq 1$. 

Let $C_Z = \cup_{z\in Z} \{z\} \times C_z$ and $D_Z = \overline{\cup_{n\geq 1}f^n C_Z}$.  Further, $f$ is vertically expanding over $Z$ if and only if $D_Z \cap J_Z = \emptyset$.  
%

As in one dimension, $z\mapsto K_z$ is upper semicontinuous and $z\mapsto J_z$ is lower semicontinuous, in the Hausdorff metric.  Further, if $f$ is vertically expanding over $Z$, then $z\mapsto J_z$ is continuous for all $z\in Z$; and if in addition, $J_z$ is connected for all $z\in Z$, then $z\mapsto K_z$ is continuous for all $z\in Z$.  However, Jonsson provides examples showing vertical expansion over $J_p \cup A_p$ neither implies that $z\mapsto J_z$ is continuous on all of $\CC$, nor that $z\mapsto K_z$ is continuous for $z\in J_p$.
We will study one example of the latter phenomenon in Section~\ref{sec:examples}.

\textbf{Axiom A polynomial skew products.}
In this paper we rely significantly on Jonsson's vertical expansion criteria for Axiom~A:

\begin{thm}[\cite{MatSk}, Theorem 8.2]
\label{thm:AxiomAVertExp}
A polynomial skew product $f$ is Axiom A on $\Ct$ if and only if 
\begin{enumerate}
\item $p$ is uniformly expanding on $J_p$,
\item $f$ is vertically expanding over $J_p$, and
\item $f$ is vertically expanding over $A_p$.
\end{enumerate}
Moreover, if $f$ is Axiom A, then $\chainrec = \Omega_f = \overline{\text{Per}(f)}$.
\end{thm} 

Note Jonsson's criticality condition for vertical expansion over $Z$ implies  
\begin{cor}[\cite{MatSk}, Corollary 8.3]
\label{cor:AxiomACritCond}
A polynomial skew product $f$ is Axiom A on $\Ct$ if and only if (1) $D_p \cap J_p =\emptyset$, (2) $D_{J_p} \cap J_2 = \emptyset$, and (3) $D_{A_p} \cap J_{A_p} = \emptyset$.
\end{cor}

Jonsson also provides a structural stability result for Axiom A skew products (\cite{MatSk}, Theorem A.6 and Proposition A.7).  A consequence is that it makes sense to refer to a connected component of the subset of Axiom A mappings in a given parameter space  as a {\em hyperbolic component}. 


\section{Examples of Axiom A polynomial skew products}
\label{sec:newclasses}


\textbf{Products.}
A simple product, $f(z,w) = (p(z),q(w))$ is Axiom A if $p$ and $q$ are hyperbolic. There are (at most) four chain transitive components: $\chainrec = (J_p \times J_q) \ \cup \ (A_p \times J_q) \ \cup \ (J_p \times A_q) \ \cup \ (A_p \times A_q)$; the first is $J_2$, the codimension zero set on which $f$ is uniformly expanding, the middle two are codimension one saddle sets, and the last set is the attracting periodic points of $f$.  For example, for $(z,w) \mapsto (z^2, w^2)$, we have $J_p \times J_q = S^1 \times S^1$ is a torus, $A_p \times J_q = \{0\} \times S^1$ and $J_p \times A_q = S^1 \times \{0\}$ are circles, and $A_p \times A_q$ is the origin.

\textbf{(Small) Perturbations of products.}
By Jonsson's structural stability results, any $f(z,w)$ which is a sufficiently small perturbation of an Axiom A product is also Axiom A, with basic sets (i.e., chain transitive components) topologically corresponding to those of the product.

\textbf{Jonsson's nonproduct example.}
In \cite{MatDil}, Diller and Jonsson describe a class of maps which are Axiom A, but not conjugate to any Axiom A product (again using Corollary~\ref{cor:AxiomACritCond}), and consider the example: 
$f_0(z,w) = (z^2 - 9900, w^2 + (z+99)/6)$.   

This type of example starts with $p(z)=z^2-R$ with $R>2$, so $J_p$ is a real Cantor set. Let $\alpha < 0 < \beta$ denote the two (necessarily repelling) fixed points of $p$.  Let $D_1 := [-\beta, -\eta]$, where $\eta = \sqrt{R-\beta}$, and let $D_2 = -D_1$. Then $J_p \subset D_1 \cup D_2$. 
For $f_0$, we have $R=9900$, $\beta=100$, $\alpha = -99$, and $\eta=\sqrt{9800} \approx 98.99$.

The skew product $f_0$ is of the form $(z^2-R, w^2 + (z+\beta)/S)$, where $S\geq 1$.  Thus for $z\in D_1$, i.e., close to $-\beta$, we have $q_z(w) \approx w^2$.  One can easily check that $\alpha+\beta=1$, so if $S$ is sufficiently large, then $q_{\alpha}(w) = w^2 + (\alpha+\beta)/S$ is conjugate to $w^2$, and so since $\alpha$ is a fixed point, $J_{\alpha}$ is a quasicircle.  On the other hand, for $z\in D_2$, i.e., $z$ close to $\beta$, we have $q_z(w) \approx w^2 +2\beta/S$, and since $\beta \geq 2$, if $S$ is not too large, then $q_z(w) \approx w^2 + c$ for some $c>0$ outside the Mandelbrot set, i.e., with escaping critical point.  Hence $J_{\beta}$ is a Cantor set.  Using the fact that $q(J_z)=J_{p(z)}$, we quickly see that the fibers $J_z$ for $z\in J_p$ contain a rich mix of Cantor sets and collections of circles.  The nonwandering set consists of $J_2$ (the closure of the repelling periodic points) and the saddle fixed point $(\alpha, 0)$.
If $R$ is sufficiently large, then such a map is Axiom A.

Diller and Jonsson describe a method for generating an Axiom A example of this type for every degree $d\geq 2$ (in fact, for every $k$ with $1 \leq k<d$, they can construct an example with topological entropy $\log k$).  $J_p$ is always contained in the union of two disjoint sets, $E$ and $E'$, with $q_z(w) \approx w^d$ for $z\in E$ and $q_z(w) = w^d + r(z)$, such that $r(z)$ is large for $z \in E'$, but $r(z)$ is of degree less than $d$.      

\textbf{Generating new examples.} 
The first proposition below is an attempt to generalize the example of Diller and Jonsson, by replacing the map $w^2$, and the circles in $J$ that arise from this map, with the map $q_c(w) = w^2+c$ for any hyperbolic $c \in \mathcal{M}$, and thus replacing circles with any connected hyperbolic Julia set $J_c$.

\begin{prop} \label{prop:matgen1}
For any hyperbolic polynomial $x \mapsto x^2+c$, with $c\in\mathcal{M}$, there exist $a\in\CC, R>6,$ and $S>0$ such that for 
\[
f(z,w) = (p(z), q_z(w)) = \(z^2-R, w^2+c + \frac{z+a}{S} \),
\]
we have:
\begin{enumerate}
\item $J_p \subset D_1 \cup D_2$, for the intervals $D_1 = [-\beta, -\eta]$ and $D_2 = -D_1$, where $\beta > 0 > \alpha$ are the fixed points of $p$,
and $\eta = \sqrt{R-\beta}$;
\item for all $z\in D_1$, $q_z$ is in the same hyperbolic component of $\mathcal{M}$ as $c$ (hence $q_z(w)$ is topologically conjugate to $w^2+c$); and
\item for all $z\in D_2$, $q_z$ is outside of $\mathcal{M}$.
\end{enumerate} 
\end{prop}

\begin{proof}
First, we quickly check (1).
The fixed points $\alpha, \beta$ are given by: $\alpha = (1-\sqrt{1+4R})/2, \beta=1+\sqrt{1+4R})/2$.
If $R>2$, then straightforward calculation gives $R > \beta > 2$, hence $\eta >0$, 
and $p^{-1} ([-\beta, \beta]) = D_1 \cup D_2$, with $p^{-1}(\beta)=\pm \beta$ and $p^{-1}(-\beta)=\pm \eta$ (and $p^{-1}(\alpha)=\pm \alpha$), and $p$ a homeomorphism on each of $D_1$ and $D_2$. The map $p(z)=z^2-R$ is conjugate to the one-sided two shift, with $J_p$ a real Cantor set contained in  $D_1\cup D_2$.
 
Let $a = (\eta + \beta)/2$. 
Now, given $a$, we find an $S$ establishing (2), for any $R>6$. Note over all the $z\in D_1$, the interval of $q_z$'s is exactly centered around $w^2+c$, with the extreme $q_z$'s differing from $w^2+c$ by length$(D_1)/(2S) = (\beta-\eta)/(2S)$.  We just need to require $S$ to be large enough to keep $q_z$ in the same hyperbolic component as $w^2+c$.
  Let $\sigma$ denote the horizontal distance between $c$ and the boundary of the hyperbolic component of $\mathcal{M}$ containing $c$.  So we need $(\beta-\eta)/(2S) < \sigma$.
Since $R>6$, $\beta-\eta=\beta - \sqrt{R-\beta} < 3-\sqrt{3}$. Hence if $S > 2\sigma /(3-\sqrt{3})$, then $(\beta-\eta)/ (2S) < \sigma$.

Finally, we show that given $a$ and $S$, we can pick $R$ large enough to satisfy (3).  For all $z\in D_2=[\eta, \beta]$, we need $q_z(w) \approx w^2 + \text{large}$, i.e., we need $q_z$ outside of $\mathcal{M}$.  The smallest $q_z$ is $q_{\eta}(w) = w^2 + c + (\eta+a)/S = w^2 + c + (3\eta + \beta)/ (2S)$.  So we just need $\Re(c) + (3\eta + \beta)/ (2S) \geq 2$ to guarantee we're outside of $\mathcal{M}$.  Since $c\in\mathcal{M}$ we know $\Re(c) \geq -2$, so we just need $(3\eta + \beta)/ (2S) \geq 4$.  But $S$ is fixed, so all we must do is make $(3\eta+\beta)$ large, which we can do by taking $R$ large enough.  Specifically, let $R$ be large enough that $\eta > 2S$. We can do this since $\eta < = \sqrt{R-\beta}$ and $\beta \approx \sqrt{R}$, so $\eta$ tends to infinity as $R$ does.
Since $\beta > \eta$, this yields the required
$(3\eta +\beta)\ (2S) > 4\eta/ (2S) = 2\eta/S > 2(2S)/S = 4$.
\end{proof}
 
\begin{rem}
We conjecture that all, or at least a large subset of, maps satisfying Proposition~\ref{prop:matgen1}  are Axiom A.
\end{rem}

Unfortunately, we do not know a general method for verifying the above remark.  We prove Axiom A for a specific example using the program {\em Hypatia} in Section~\ref{sec:examples}. 

One would ideally generalize this class of examples further, by replacing the Cantor sets above with an arbitrary hyperbolic polynomial, i.e., construct a skew product of the form $(z^2-R, w^2+ l(z))$, with $l(z)$ some linear or quadratic map, but instead of (3) above require $q_z$ for $z\in D_2$ to be in the same hyperbolic component of $\mathcal{M}$ as some other hyperbolic map $x^2+c_2$.  

However, this does not work in general for any two hyperbolic maps $x^2+c_1, x^2+c_2$.  
For, we would need $l(z)$ to map $[-\beta, \beta]$ into the parameter space of $c$'s, so that $l(D_i)$ is contained in the hyperbolic component of $\mathcal{M}$ containing $c_i$.  But $R>2$ is necessary, and so the length of $D_i$ is not too small, but the hyperbolic components containing $c_i$ may be small compared to the distance between $c_1$ and $c_2$.  Hence there may be no linear or quadratic map $l(z)$ which will work.  

Hence the best we can do is consider $c_1, c_2 \in \mathcal{M}$ hyperbolic, ideally near the center of their hyperbolic components, then choose some $R>6$ and consider 
\begin{equation} \label{eqn:gencls2}
f(z,w) = (p(z), q_z(w)) = \(z^2-R, \ \ w^2+ \frac{c_1}{2}\(1-\frac{z}{a}\) + \frac{c_2}{2}\(1+\frac{z}{a}\) \),
\end{equation}
where $a = (\eta + \beta)/2$ is the center of $D_2$ as in Proposition~\ref{prop:matgen1}.  
Then $q_z(w)$ is of the form $w^2 + l(z)$, and $l(z)$ is the linear map satisfying $l(-a) = c_1$ and $l(a)=c_2$.  Then we must check whether $l(-\beta)$ and $l(-\eta)$ are in the same hyperbolic component as $c_1$, and whether $l(\eta)$ and $l(\beta)$ are in the same hyperbolic component as $c_2$.  If this holds, the $l(D_i)$ is in the same hyperbolic component as $c_i$, and we have the desired conclusion.  In this case, we conjecture the map is Axiom A.

Fortunately, for the several largest hyperbolic components of $\mathcal{M}$, the above construction is successful. For example:

\begin{exmp}[Circles and Basillicas]
\label{exmp:circxbas}
For $f(z,w) = (p(z), q_z(w))= (z^2-90, w^2 - z/20 - 1/2)$, it is easy to check that:
(1) $J_p \subset D_1 \cup D_2$, with $D_1=[-10,-\sqrt{80}]$;
(2) for all $z\in D_1$, $q_z(w)$ is in the same hyperbolic component as $w^2$; and
(3) for all $z\in D_2$, $q_z(w)$ is in the same hyperbolic component as $w^2-1$ (the Basillica).
\end{exmp}

\begin{exmp}[Basillicas and Rabbits]
\label{exmp:basxrab}
For $f(z,w) = (p(z), q_z(w))= (z^2-90, w^2 + z(0.046 + 0.04i) -.56+.37i)$, it is easy to check that:
(1) $J_p \subset D_1 \cup D_2$, with $D_1=[-10,-\sqrt{80}]$;
(2) for all $z\in D_1$, $q_z(w)$ is in the same hyperbolic component as $w^2-1$ (the Basillica); and
(3) for all $z\in D_2$, $q_z(w)$ is in the same hyperbolic component as $w^2-0.12+0.75i$ (the Rabbit).
\end{exmp}

\section{{\em Hypatia} for polynomial skew products}
\label{sec:hypatiatheory}

In this section we describe the algorithms for studying polynomial skew products that we have
 implemented in {\em Hypatia}.  There are two main steps in {\em Hypatia}: 
 \begin{enumerate}
\item building a neighborhood $\boxcov$ of some $\chainrec_{\delta}$ and a graph modelling the 1-step dynamics of $f$ on $\boxcov$,
and
\item testing whether a skew product $f$ is Axiom A, using Jonsson's vertical expansion criteria, by attempting to build a metric for which $f$ is vertically expanding over neighborhoods of $J_p$ and $A_p$.
\end{enumerate}

\subsection{Building a model of the chain recurrent set}
\label{sec:chrechypatia}

In \cite{SLHone} we describe an algorithm we call the {\em box chain construction}, 
for building a directed graph $\Gamma$ modelling a given map $g$ of $\Cn$ on a 
neighborhood $\boxcov$ of $\chainrec$.  
A similar approach in different settings can be found in \cite{DellHoh, Dell1, Eiden, KMisch,
Osi, OsiCamp}. The basic output of this construction is a graph called a {\em box chain recurrent
model}, which we define below.

\begin{defn} \label{defn:boxchrecmodel}
Let $\Lambda$ be an invariant set of a map $g \colon \Cn \to \Cn$.
Let $\Gamma = (\mathcal{V}, \mathcal{E})$ be a directed graph,
with vertex set $\mathcal{V} = \{ B_k \}_{k=1}^N$,  a finite collection of closed boxes in $\Cn$, having disjoint interiors, and  such that the union of the boxes $\boxcov = \cup_{k=1}^N B_k$ contains $\Lambda$.  Suppose there is a $\delta>0$ such that 
$\Gamma$ contains an edge from $B_k$ to $B_j$ if the image $g(B_k)$ intersects a $\delta$-neighborhood of $B_j$, \textit{i.e.},  
\[
\mathcal{E} \supset \{ (k, j) \colon g(B_k) \cap \mathcal{N} (B_j,{\delta} ) \neq \emptyset \}.
\]
Further, assume $\Gamma$ is partitioned by edge-connected components $\Gamma'$ which are strongly connected, i.e., for each pair of vertices $B_k, B_j$, there is a path in $\Gamma'$ from $B_k$ to $B_j$, and vice-versa.   Then we say $\Gamma$ is a {\em box chain recurrent model} of $g$ on $\Lambda$, and each $\Gamma'$ is a {\em box chain transitive component}.
\end{defn}

A \boxchrecmod $\Gamma$ is an approximation to the dynamics of $g$ on $\chainrec$, 
and the connected components $\Gamma'$ of $\Gamma$ are approximations to the chain transitive components.  

The basic box chain construction for a polynomial map $p$ of $\CC$ can be summarized as follows:
\begin{enumerate}
\item Compute $R$ such that for some $\delta>0$, $\chainrec_{\delta} \subset \boxcov_0 := [-R,R]^2$.
\item Subdivide $\boxcov_0$ into a grid of boxes $\{ B_k\}$.
\item Build a graph $\Upsilon$ with vertices $\{ B_k \}$ and edges $\{ (k,j) \colon p(B_k) \cap \mathcal{N} (B_j,{\delta} ) \neq \emptyset \}$. 
\item Find the maximal subgraph $\Gamma$ of $\Upsilon$ which consists precisely of  edges and vertices lying in cycles.  Then $\Gamma$ is partitioned by its edge-connected components: $\Gamma = ( \Gamma^1 \sqcup \Gamma^2 \sqcup \ldots \sqcup \Gamma^N )$.
\item If desired, refine $\Gamma$ by subdividing the boxes of $\Gamma$ and repeating (3) and (4).
\end{enumerate}
Then $\Gamma$ is a \boxchrecmod of $p$ on $\chainrec$, and the $\Gamma^j$ are the box chain transitive components.

Now we outline the process we follow for skew products.  

(1) Use the one dimensional box chain construction of \cite{SLHone} to build a box chain recurrent model $\Gamma$ for $p$ in the $z$-plane.  Iterate the construction until the boxes are small enough that (a) $\Gamma$ separates the box chain transitive component containing $J_p$, call it $\Gamma^0_{J_p}$, from the one containing $A_p$ (the attracting periodic orbit, if there is one), called $\Gamma^0_{A_p}$ and (b) hyperbolicity of $p$ can be established 
for $\Gamma^0_{J_p}$ (using the techniques of \cite{SLHtwo}).  Construction begins by calculating an $R_1$ such that $\chainrec_p \subset [-R_1, R_1]^2$.  The boxes of $\Gamma$, $\{B^z_k\}$, are a subset of the boxes in a $2^n \times 2^n$ grid on $[-R_1, R_1]^2$.

(2) Build a model of $f$ in only the box fibers over $\Gamma^0_{A_p}$ and $\Gamma^0_{J_p}$.  

That is, first compute an $R_2$ so that $\mathcal{R}_{\delta}(f) \subset [-R_1, R_1]^2 \times [-R_2, R_2]^2$ (this calculation is similar to computing $R_1$, which is explained in \cite{SLHone}).
Then choose some $m>0$ and construct a grid of $2^m \times 2^m$ boxes $\{B^w_j\}_{j=1}^{2^{2m}}$ in the square $[-R_2, R_2]^2$ in the $w$-plane.  
Form boxes in $\Ct$ as the product of boxes in $\Gamma^0_{A_p}$ or $\Gamma^0_{J_p}$ times boxes in the new $w$-grid, i.e.,
for each box $B^z_k$ in $\Gamma^0_{A_p}$ or $\Gamma^0_{J_p}$,  
we have $2^m \times 2^m$ boxes in $\Ct$ of the form $B_{k,j} = B^z_k \times B^w_j$. Thus we get a set of boxes in $\Ct$, $\{ B_l \}$.  

Now build a transition graph for $f$, $\Upsilon$, whose vertices are the boxes $B_l$ and with an edge $B_k \to B_j$ if $f(B_k) \cap \mathcal{N} (B_j,{\delta} ) \neq \emptyset$.  Finally, find the subgraph $\Gamma$ of $\Upsilon$ consisting precisely of the vertices and edges of $\Upsilon$ which lie in cycles, and decompose $\Gamma$ into its edge-connected components: $\Gamma = ( \Gamma^1 \sqcup \Gamma^2 \sqcup \ldots \sqcup \Gamma^N )$.  Each $\Gamma^j$ is a box chain transitive component.  
If the boxes in $\Ct$ are sufficiently small, then the box chain transitive components separate the chain transitive components of $f$; for instance, any box chain transitive components over $\Gamma^0_{J_p}$ should be distinct from those over $\Gamma^0_{A_p}$.    There may be more than one box chain transitive component over an invariant set in the base, for example if $q(z,w)=w^2$, then there is a component containing $S^1=J_z$ in the $w$-plane, and (if boxes are small enough) a separate component containing the attracting fixed point of $q_z$, the $w$-plane's origin.

Recall $J_Z = \overline{\cup_{z\in Z}\{z\} \times J_z }$.
Let $\Gamma^1_{A_p}$ denote the  box chain transitive component over $\Gamma^0_{A_p}$ which contains $J_{A_p}$; similarly let $\Gamma^2_{J_p}$ be the box chain transitive component over $\Gamma^0_{J_p}$ containing $J_{J_p}$.   Heuristically, we can quickly identify the box chain transitive component over a set $Z$ which contains $J_Z$, because it will be the largest by far (i.e., have the most vertices and edges).  

\subsection{Establishing Axiom A: box vertical expansion}
\label{sec:showaxiomA}

After the modified box chain construction of the previous section has produced box chain transitive components  $\Gamma^1_{A_p}$ (containing $J_{A_p}$) and $\Gamma^2_{J_p}$ (containing $J_{J_p}$), it is time to test whether $f$ is Axiom A.  We use the vertical expansion criteria of Theorem~\ref{thm:AxiomAVertExp}, i.e., we want (1) $p$ uniformly expanding on $J_p$, (2) $f$ vertically expanding over $J_p$, and (3) $f$ vertically expanding over $A_p$.  Note (1) was checked in the construction of the previous subsection, so in this subsection we deal with (2) and (3).

Recall that $f$ is vertically expanding over $Z$ if there are $c>0, \lambda >1$ such that $\abs{DQ^n_z(w)} \geq c\lambda^n$, for all $z\in Z$, $w\in J_z$, and $n\geq 1$, where let  $z_n = p^n(z)$, $q_z(w) = q(z,w)$, and $Q^n_z(w) = q_{z_{n-1}} \circ \cdots \circ q_z$, so that $f^n(z,w) = (z_n, Q^n_z(w)).$   Thus vertical expansion is expansion in one complex dimension only, specifically, in the $w$-coordinate.
This allows us to easily adapt our algorithm from \cite{SLHtwo} for verifying hyperbolicity of polynomial maps of $p$ to testing vertical expansion.

In one dimension, hyperbolicity of $p$ reduces to uniform expansion of tangent vectors over a neighborhood of $J$, in some riemannian metric.  {\em Hypatia} tests hyperbolicity by trying to build a piecewise constant metric in which $p$ is expanding:

\begin{defn} \label{defn:boxexp}
Let $p$ be a polynomial map of $\CC$.
Let $\Gamma'$ be the box chain transitive component containing $J_p$, with box vertices $\{B_k\}_{k=1}^N $.   
Call $p$ {\em box-expansive} on $\Gamma'$
if there exists $L>1$ and positive constants $\{\phi_k \}_{k=1}^N$ such that for every edge $(k,j)$ in $\Gamma'$, we have $\phi_j \abs{D_z p(v)} \geq L \phi_k \abs{v}$, for all $z\in B_k$, and all $v\in T_z\CC$ (or simply $\phi_j \abs{p'(z)} \geq L \phi_k$).
\end{defn}

The constants  $\{\phi_k \}_{k=1}^N$ define what we call a {\em box metric} on $\boxcov' = \cup_{k=1}^N B_k$, i.e., simply a piecewise constant multiple of the euclidean metric on $\boxcov'$, with a different constant for each box. In \cite{SLHtwo}, we establish the following theorem, showing that box expansion implies the standard notion of hyperbolicity:

\begin{thm}\label{thm:boxexp}
Suppose $p$ is box-expansive by $L>1$ on $\boxcov' = \cup_{k=1}^N B_k$.  Let $\Lambda_1 = \boxcov' \cap p^{-1}(\boxcov')$. Then there exists an $\eta>1$ and a continuous norm $\abs{\cdot}_{\rho}$ on $T_{\Lambda_1}\CC$ which $p$ expands by $\eta$, i.e., $\abs{D_zp(v)}_{\rho} \geq \eta \abs{v}_{\rho}$, for all $z\in\Lambda_1$, and all $v \in T_z\CC$.
\end{thm}

The notion of box expansion is easily adapted to vertical expansion.

\begin{defn}\label{defn:boxvertexp}
Let $f(z,w) = (p(z), q(z,w))$ be a polynomial skew product.
Let $Z \subset K_p$ be compact with $pZ \subset Z$.
Let $\Gamma'$ be the box chain transitive component containing $J_Z = \overline{\cup_{z\in Z}\{z\} \times J_z }$, with box vertices $\{B_k\}_{k=1}^N $.   Call $f$ {\em box vertically expansive} on $\Gamma'$ if there exists $L >1$ and positive constants $\{\phi_k \}_{k=1}^N$ such that for every edge $(k,j)$ in $\Gamma'$, we have $\phi_j \abs{D_w q_z(v)} \geq L \phi_k \abs{v}$, for all $(z,w) \in B_k$, and all $v \in T_w \CC$ (or simply $\phi_j \abs{q'_z(w)} \geq L \phi_k$).
\end{defn}

Theorem~\ref{thm:boxexp} immediately yields as a corollary its analog for vertical expansion. 
We use that the definition of vertical expansion over $Z$ given by Jonsson is equivalent to the existence of some constant $\eta>1$ and some riemannian norm $\abs{\cdot}_{\rho}$ in the tangent bundle to the $w$-coordinate over $Z$, such that in this norm, $DQ_z$ expands tangent vectors by $\eta$.
That is, if $J_Z \subset 
\boxcov' = \cup_{B_k \in \Gamma'} B_k$, then box vertical expansion for $\Gamma'$ implies vertical expansion of $f$ over a neighborhood of $Z$, in some (continuous) riemannian norm, on the $w$-coordinate tangent bundle for all $(z,w)$ in a neighborhood $\Lambda_1$, with $J_Z \subset \Lambda_1 \subset \boxcov'$.

We describe an algorithm for attempting to establish box expansion in \cite{SLHtwo}.  
To test box vertical expansion, we can use this algorithm with only very small changes.
We perform the one-dimensional algorithm twice, once for $\Gamma'  = \Gamma^1_{A_p}$ (containing $J_{A_p}$) and then again for $\Gamma' = \Gamma^2_{J_p}$ (containing $J_{J_p}$).  For each $\Gamma'$, we use the one-dimensional recursive algorithm to attempt to build a set of metric constants $\{ \phi_k \}$ for each box $B_k$ of $\Gamma$, such that for each edge $(k,j)$ in $\Gamma'$, the inequality 
\[ 
\phi_j \geq \frac{\phi_k \ L}{\displaystyle\inf_{(z,w)\in B_k} \abs{q'_z(w)}}
\]
 is satisfied. 
Note if any of the derivatives $q_z'(w)$ is zero, then $\Gamma'$ contains a critical point (of $q_z$), so box vertical expansion necessarily fails and the algorithm is terminated.  At this point, the user may re-try with smaller boxes, to attempt to separate this critical point from $J_Z$.  If all such derivatives are positive, the recursive algorithm to define the $\phi_k$ may still fail.  This will happen if  the map is not vertically expanding over $Z$, or it may mean the boxes or $L$ are simply too large.  Thus the user may re-try with smaller L, or smaller boxes.  
On the otheer hand, if the algorithm is successfully completed, producing a set of metric constants $\{ \phi_k \}$ for $\Gamma'$, then Theorem~\ref{thm:boxexp} implies that $f$ is vertically expanding over $J_p$ and $A_p$, i.e., (2) and (3) of Theorem~\ref{thm:AxiomAVertExp} are satisfied, and combined with (1), which we checked in Section~\ref{sec:chrechypatia}, this yields that $f$ is Axiom A.

\section{Results from {\em Hypatia}}
\label{sec:examples}

In this section, we describe the results of applying the algorithms of Section~\ref{sec:hypatiatheory} to some specific examples of maps from the classes discussed in Section~\ref{sec:newclasses}.

The computations described in this section were run via a C$^{++}$ program in a Unix environment,
with 2GB of RAM.  When computations became
overwhelming, memory usage was the limiting factor.  
\begin{footnote}{
All computations were performed on the Indiana University IBM Research SP system.  Hence, this work was supported in part by Shared University Research grants from IBM, Inc. to Indiana University.
}\end{footnote}
Rigor was maintained in the computations using the {\em Interval Arithmetic} routines in the PROFIL/BIAS package, available at \cite{PBIA}.  See Appendix~\ref{sec:IA} for a brief introduction to this technique.

The examples we examine are all quadratic skew products.  Thus $p$ and $q_z$ each have only one critical point, at the origin, and $p$ has at most one attracting cycle.   

Before proceeding with the examples, we summarize the data given.  For each example, the box chain transitive components were formed by boxes from a $(2^n)^2 \times (2^m)^2$ grid on $\boxcov = [-R_1, R_1]^2 \times [-R_2, R_2]^2$.  For each graph component for which vertical expansion was tested, we give the size of the graph, $(V; E)$,  with $V$ the number of vertices and $E$ the number of edges. If expansion was established, we state the expansion constant, $L$, the range of metric constants $\phi_k$ constructed, as an interval $[\phi_{\min}, \phi_{\max}]$, and the average metric constant $\phi_{\text{avg}}$.

\begin{exmp}[Circle cross circle perturbation]
\label{exmp:circxcirc}
The map $f(z,w) = (z^2, w^2 + w/10 + z/100)$
 appears to be a stable perturbation of the product $(z,w) \mapsto (z^2, w^2)$.  
 This map is in Heinemann's family of {\em cannelloni} maps, since $J_p$ is connected and the fiber Julia sets are quasicircles.
 
 We proved with {\em Hypatia} that it is Axiom A, using boxes from a $(2^5)^2 \times (2^4)^2$ grid on $\boxcov = [-1.1, 1.1]^2 \times [-1.21, 1.21]^2$.
Since $z \mapsto z^2$ has an attracting fixed point at the origin, there are two box chain transitive components to test expansion over: one for $J_{J_p}$ and one for $J_{A_p} = \{ 0 \} \times J_{ \{0\} }$.
The component containing $J_{J_p}$, of size $(V;E)=(39{,}200; 1{,}953{,}900)$,
is box vertically expansive for $L=1.29105$,
with metric constants in $[0.1674,1]$, with average $0.290$.
The component containing $J_{A_p}$, of size $(V;E)=(608; 4{,}416)$, is box vertically expansive for $L=1.29105$, with metric constants in $[0.11116,1]$, with average $0.168$.
%


\end{exmp}

\begin{exmp}[Cantor cross circle perturbation]
\label{exmp:cantxcirc}
For the map $f(z,w) = (z^2 + 2, w^2 + z/10)$, we have $J_p$ is a Cantor
set in the base, and $q_z$ is a small perturbation of $w \mapsto w^2$.
Heinemann calls this type of perturbation a {\em Cantor skew},  and shows for sufficiently small perturbation, the fiber Julia sets are Jordan curves.  Jonsson mentions they are in fact quasicircles, and shows, using Corollary~\ref{cor:AxiomACritCond}, that Cantor skews are Axiom A.
Note such a map has only one basic set, which is topologically a Cantor set crossed with a circle.

For this map, we established Axiom A with boxes from a $(2^7)^2 \times (2^5)^2$ grid on $\boxcov = [-2.1, 2.1]^2 \times [-1.28, 1.28]^2$.  The
 box chain transitive component containing $J_{J_p}$ was of size $(V;E)=(39{,}880; 1{,}447{,}640)$.
We established box vertical expansion on this component by $L=1.125$, with metric constants
in $[0.0024,1]$, with average $0.00529$.  
See Figure~\ref{fig:cantxcirc}.

\begin{figure}
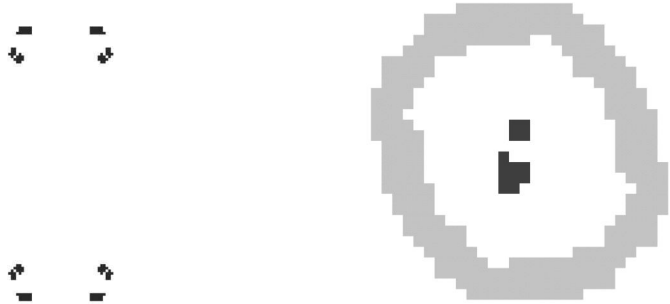

\begin{center}
\drawfigzcmplxcant \hskip3cm \drawfigcantcirc 
\caption{Let $f(z,w) = (z^2 + 2, w^2 + z/10)$, as in Example~\ref{exmp:cantxcirc}. On the left is the box cover of $J_p$, with boxes from a $(2^7)^2$ grid on $[-2.1, 2.1]^2$.  On the right is the slice of the box chain transitive components in the fiber over a repelling fixed point of $p$,  $z=0.5 + 1.32288i$ (necessarily in $J_p$), with boxes in the $w$-coordinate from a $(2^5)^2$ grid on $[-1.28, 1.28]^2$. Since $q_z(w)$ is close to $w^2$, there are two components: one containing the quasicircle $J_z$, and the other a collection of boxes near the origin.} 
\label{fig:cantxcirc}
\end{center}
\end{figure}
\end{exmp}

\begin{exmp}[Cantor cross Basillica perturbation]
\label{exmp:cantxbas}
One can generalize the above example by examining $J_p$ Cantor, and $q_z$ a small perturbation of any connected, hyperbolic polynomial, like the Basillica: $w \mapsto w^2-1$.

For the specific example $f(z,w) = (z^2-6, w^2-1-z/100)$, we achieved separation of the chain transitive components using boxes from a $(2^{10})^2 \times (2^7)^2$ grid on $\boxcov = [-3.1, 3.1]^2 \times [-1.732, 1.732]^2$ (see Figure~\ref{fig:cantxbas}), but were unable to show Axiom A at that level, or even for a refinement of boxes from a  $(2^{10})^2\times (2^8)^2$ grid on $\boxcov$. Trying to refine further maxed out our memory resources of 2GB of RAM.

\begin{figure}
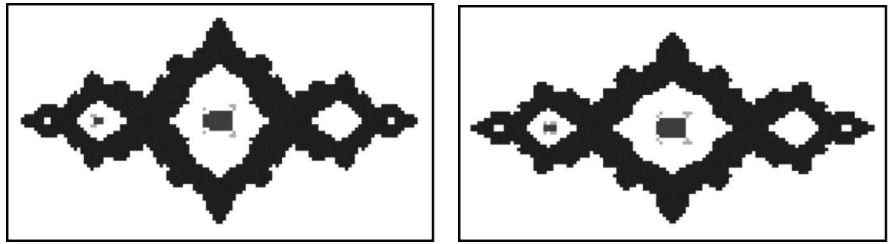

\begin{center}
\drawfigcantxbasnegb \ \ \drawfigcantxbasb 
\caption{For the map $f(z,w) = (z^2-6, w^2-1-z/100)$ of Example~\ref{exmp:cantxbas}, $J_p$ is a Cantor set, while each $q_z$ is close to a Basillica.  Shown above are slices of the box chain transitive components in two  fibers:  left over $z=-3=-\beta$, right over $z=3=\beta$. Both show the same window, with boxes from a  $(2^7)^2$ grid on $[-1.732, 1.732]^2$ in the $w$-plane.  Different box chain transitive components are colored different shades.  There is one Basillica-like component containing $J_{J_p}$, a second component near the  attracting 2-cycle of the Basillica, $\{ 0, -1\}$, and a third component of islands skipping around the second component.  This last component contains {\em pseudo}-recurrent points only, and disappears for smaller boxes.}
\label{fig:cantxbas}
\end{center}
\end{figure}
\end{exmp}

\begin{exmp}[Jonsson's nonproduct] \label{exmp:matcirc}
The map $f(z,w) = (z^2-90, w^2+ (z+9)/4)$ is an example of the class of map which  Jonsson showed to be Axiom A, but not conjugate to any product.  $J_z$ in the fiber over the $\alpha$-fixed point of $p$ is a circle, while $J_z$ over the $\beta$-fixed point of $p$ is a Cantor set.

For this specific map, we used {\em Hypatia} to prove Axiom A with
boxes from a $(2^7)^2 \times (2^7)^2$ grid on $\boxcov=[-10.1, 10.1]^2 \times [-2.842, 2.842]^2$.  
Since $J_p$ is a Cantor set, there is only one box chain transitive component for which we must establish vertical expansion: the one containing $J_{J_p}$.
For a component of size $(V;E)=(5{,}640; 209{,}572)$, 
we showed box vertical expansion by $L=1.25$, with metric constants in $[0.066, 1]$, with average $0.125$.
See Figure~\ref{fig:matcirc}.

\begin{figure}
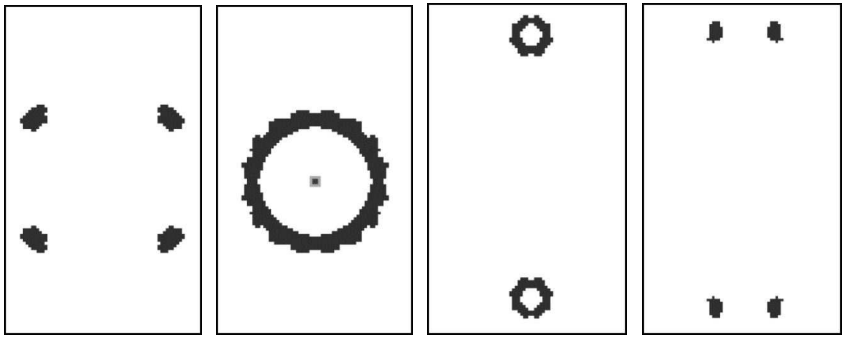

\begin{center}
\drawfigmatcircnegb \hskip.2cm \drawfigmatcirca \hskip.2cm \drawfigmatcircnega \hskip.2cm \drawfigmatcircb
\caption{For Example~\ref{exmp:matcirc}, the map $(z,w) \to (z^2-90, w^2+ (z+9)/4)$, shown are 
box chain transitive components in four fibers: left to right over $z=-10=-\beta, z=-9=\alpha, z=9=-\alpha,$ and $z=10=\beta$, with boxes from a $(2^{7})^2$ $w$-grid on $[-2.842, 2.842]^2$.  The same window is shown for each fiber.}
\label{fig:matcirc}
\end{center}
\end{figure}
\end{exmp}

\begin{exmp}[Generalization type 1 of Jonsson's nonproduct] \label{exmp:matrab}
The map $f(z,w) = (z^2-90, w^2 + z/6 + 1.4 + 0.75i)$ satisfies Proposition~\ref{prop:matgen1}.  For this map, $J_z$ in the fiber over the $\alpha$-fixed point is the Rabbit, while $J_z$ in the fiber over the $\beta$-fixed point is a Cantor set.   Note since $A_p$ is empty, we need only test vertical expansion on the  box chain transitive component containing $J_{J_p}$.
For this map, we used {\em Hypatia} to prove Axiom A, with
boxes from a $(2^{10})^2 \times (2^9)^2$ grid on $\boxcov = [-10.1, 10.1]^2 \times [-2.426, 2.426]^2$. 

See Figure~\ref{fig:matrab}.  Note that the attracting 3-cycle of the Rabbit yields a saddle $3$-cycle in $\Ct$, lying in the fiber over the $\alpha$ fixed point .  This 3-cycle is contained in its own box chain transitive component , which must be separated from the component containing $J_{J_p}$ before vertical expansion can be established.   The component of $J_{J_p}$, of size $(V;E)=(78{,}994; 2{,}066{,}558)$, is box vertically expansive by $L=1.0625$, with metric constants in $[0.000853, 1]$, with average $0.00315$.

\begin{figure}
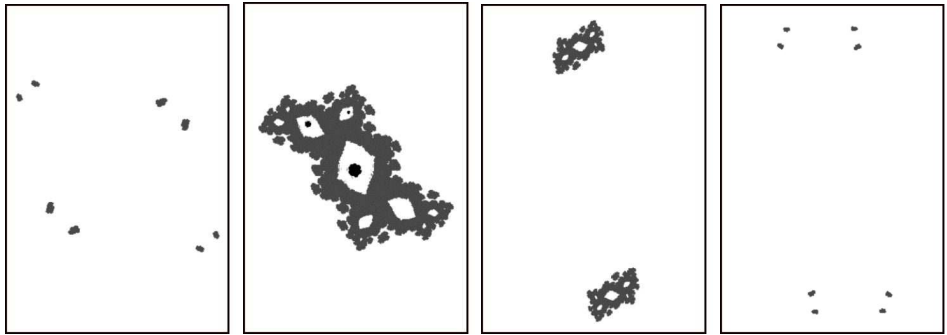

\begin{center}
\drawfigmatrabnegb \hskip.1cm \drawfigmatraba \hskip.1cm \drawfigmatrabnega \hskip.1cm \drawfigmatrabb
\caption{For Example~\ref{exmp:matrab}, the map $(z,w) \to (z^2-90, w^2+ z/6 + 1.4 + 0.75i)$, shown are 
box chain transitive components in four fibers: left to right over $z=-10=-\beta, z=-9=\alpha, z=9=-\alpha,$ and $z=10=\beta$, with boxes from a $(2^{9})^2$ $w$-grid on $[-2.426, 2.426]^2$.  The same window is shown for each fiber.}
\label{fig:matrab}
\end{center}
\end{figure}

\end{exmp}

\begin{exmp}[Generalization type 2 of Jonsson's nonproduct]
Example~\ref{exmp:circxbas} is the simplest map of this type, with $J_p$ a real Cantor set,
and the fibers containing (among other things) circles and Basillicas.
We tested {\em Hypatia} on the map $f(z,w)=(z^2-90, w^2-z/20-1/2)$ of this type.  However, we were unable to achieve separation of all of the chain transitive components, hence could not test for Axiom A.  The obstacle was memory usage.  The deepest level we could completely build  used boxes from a $(2^{11})^2 \times (2^{8})^2$ grid on $[-10.1, 10.1]^2 \times [-1.7203,1.7203]^2$, but for these boxes the critical point $0$ of $q_z$ was still in the box chain transitive component for $J_{J_p}$, and we used our entire 2GB available trying to build the next level.
\end{exmp}

\subsection{Conclusions and remaining questions}
\label{sec:conc}

In conclusion, we would like to summarize our results, and provide the reader with questions
for further study.  We have begun to examine what kind of dynamics can 
occur for Axiom A skew products, by
developing and implementing rigorous computer algorithms into a program {\em Hypatia}, to test whether a skew product of the form $(z^2+a, w^2+bw+cz+e)$ is Axiom A.
We have used {\em Hypatia} to verify the Axiom A property for Heinemann's and Jonsson's Examples; and further, derived formulas for two new classes of proposed Axiom A maps (generalizing Jonsson's examples), and used {\em Hypatia} to prove Axiom A for some examples in our first new class.

Polynomial skew products are a rich area for study.  We would like to encourage the interested reader to begin an exploration of this family, by providing the following list of questions
which arose during our recent introduction to this area. 

\begin{prob}
Produce a combinatorial description of the dynamics of maps satisfying Proposition~\ref{prop:matgen1}, or of the same type as Examples~\ref{exmp:circxbas} and~\ref{exmp:basxrab}.  For example, by
employing the theory of ``random'' dynamical systems, and studying a random composition of two hyperbolic polynomials of $\CC$ (like $z^2$ and $z^2-1$ for Example~\ref{exmp:circxbas}).  
\end{prob}

\begin{que} \label{que:fiberhomeo}
If a polynomial skew product is an Axiom A stable perturbation of a product, then are all of the fibers homeomorphic? 
\end{que}

\begin{que}
Is there an Axiom A polynomial skew product with $J_p$ connected
but whose fiber Julia sets $J_z$ for $z\in J_p$ are not all homeomorphic?  Jonsson (\cite{MatSk}, Example 9.6) and Proposition~\ref{prop:matgen1} both provide examples
of Axiom A maps whose fibers are not homeomorphic, but $J_p$ is a Cantor set.
\end{que}

\begin{que}
What can be said about polynomial skew products using the techniques of holomorphic motions, in dynamical space and/or parameter space?
\end{que}

\begin{que}
What further measure theoretic results hold for skew products?  For example, is there a ``typical'' $J_z$, for $z\in J_p$, in the measure theoretic sense?
\end{que}

\begin{que}
How can the theory of partial hyperbolicity, which is well developed in the real variables setting, and can pertain to skew products of real manifolds, shed light on the structure of polynomial skew products of $\Ct$?
\end{que}

\begin{que}
Jonsson (\cite{MatSk}, Theorem A.6) establishes structural stability  on $J_{J_p}$,
 for Axiom A polynomial skew products.  When (if ever) can structural stability be
extended to hold over a larger set?
\end{que}

\appendix

\section{Rigorous Arithmetic}
\label{sec:IA}

On a computer, we cannot work with real numbers; instead we work over
the finite space $\FF$ of numbers representable by binary floating point
numbers no longer than a certain length.  
For example, since the number 
$0.1$ is not a dyadic rational, it has an infinite binary expansion.  Thus 
the computer cannot encode $0.1$ exactly. 
Interval arithmetic (IA) provides a method for maintaining rigor in computations, and also is natural and efficient for manipulating boxes.  The basic objects of IA
 are closed intervals, $[a] = [\underline{a}, \bar{a}] \in \II\KK$, with end points in some fixed
 field, $\KK$.   An arithmetical operation on two intervals produces a resulting interval
which contains the real answer.  For example,
\[
[a] + [b] := \left[ \underline{a}+\underline{b}, \bar{a}+\bar{b} \right]
\hspace{0.5cm} \text{and} \hspace{0.5cm}
[a] - [b] :=  \left[ \underline{a}-\bar{b}, \bar{a}-\underline{b} \right] 
\]
Multiplication and division can also be defined in IA.

Since an arithmetical operation on two computer numbers in $\FF$ may not have a 
result in $\FF$,  in order to implement rigorous IA we must round
outward the result of any interval
arithmetic operation, \textit{e.g.} for $[a],[b] \in \II\FF$,
$$
[a] + [b]  := 
\left[ \left\downarrow \underline{a}+\underline{b} \right\downarrow, 
         \left\uparrow  \bar{a}+\bar{b} \right\uparrow \right],
$$
where $\left\downarrow x \right\downarrow$ denotes the largest number in $\FF$ 
that is 
strictly less than $x$ (\textit{i.e.}, $x$ rounded down), and 
 $\left\uparrow x 
\right\uparrow$ 
denotes the smallest number in $\FF$ that is strictly greater than $x$ 
(\textit{i.e.}, $x$ 
rounded up).  This is called IA with \textit{directed rounding}.

For any $x\in \RR$, let Hull$(x)$ be the smallest interval in $\FF$ which contains $x$.
That is, if $x \in \FF$, then Hull$(x)$ denotes $[x, x]$.  If $x \in \RR \setminus \FF$, then 
 Hull$(x)$ denotes $\left[ \left\downarrow x \right\downarrow, 
\left\uparrow x \right\uparrow \right]$.
Similarly, for a set $S\subset \RR$, we say Hull$(S)$ for the smallest  interval containg $S$. Whether Hull$(S)$ is in $\IR$ or $\IF$ should be clear from context.

In higher dimensions, IA operations can be carried out component-wise, on 
{\em interval vectors}.  
So if $x\in \Rn$, then Hull$(x) = \text{ Hull}(x_1) \times \cdots \times
\text{ Hull}(x_n)$, and if $S \subset \Rn$, then $\text{Hull}(S)$ is the smallest vector in $\IF^n$ (or $\IR^n$) containing $S$.
Note to deal with intervals in $\Cn$ we simply identify $\Cn$ with $\RR^{2n}$. 
Thus a box in $\CC = \Rt$ is an interval vector of length two.
Our extensive use of boxes is designed to 
make IA calculations natural.

To compute the image of a point $x$ under a map $f$ using IA, first convert $x$ to an interval vector $X = \text{ Hull}(x)$, then use a (carefully chosen) combination of the basic arithmetical operations to compute an interval vector $F(X)$, such that we are guaranteed that Hull$(f(x)) \subset F(X)$.  
If $f \colon \Rn \to \Rn$ is continuous, then an \textit{interval extension} of $f$, $F = F(f)$, is a function which maps a box $B$ in $\Rn$ to a box $F(B)$ containing $f(B)$, \textit{i.e.},  $F(B)\supset \text{Hull}(f(B))$.
Usually, 
we would like $F(B)$ to be as close as possible to Hull($f(B)$).
We shall not discuss here how to find the best $F$. 

Each time an arithmetical calculation is performed, one must think carefully
about how to use IA. For example, IA is not distributive.  Also, it can easily create large error propagation.
For example, iterating a polynomial map on an interval vector (which is not very close to an 
attracting period cycle) 
will produce a very large interval vector after only a few iterates.  
That is, if 
$B$ is a box in $\Rt = \CC$, and one attempts to 
compute a box containing $p_c^{10}(B)$, for $p_c(z)=z^2+c$, by:
\begin{tabbing}
\hspace*{.25in} \= \hspace{.25in} \= \hspace{.25in} \kill
\>for $j$ from $1$ to $10$ do \\
\>\>$B = p_c(B)$
\end{tabbing}
then the box $B$ will likely grow so large that its defining bounds become 
machine $\infty$, \textit{i.e.}, the largest floating point in $\FF$.

For further background on interval arithmetic, see  \cite{ GenIA, MooreIA1, MooreIA2}.
We use IA for all the rigorous computations in the computer program {\em Hypatia}.  The
IA routines were all provided by the PROFIL/BIAS package, available at 
\cite{PBIA}.  


\bibliographystyle{plain}
\bibliography{Skew}

\end{document}